\documentclass[12pt,a4paper]{amsart}
\usepackage[utf8]{inputenc}
\usepackage{amsmath}
\usepackage{amsfonts}
\usepackage{amssymb}

\newtheorem{theorem}{Theorem}           
\newtheorem{lemma}{Lemma}               
\newtheorem{assumption}{Assumption}
\newtheorem{definition}{Definition}
\newtheorem{remark}{Remark}

\newcommand{\RR}{\mathbb R}
\newcommand{\NN}{\mathbb N}
\newcommand{\ZZ}{\mathbb Z}

\begin{document}

\title[On maximally oscillating perfect splines] {On maximally oscillating perfect splines and some of their extremal properties}

\author{Oleg Kovalenko}

\address{Oleg Kovalenko, Department of Mathematics and Mechanics, Oles Honchar Dnipro National University, Gagarina~ave., 72, Dnipro, 49010, Ukraine}

\email{olegkovalenko90@gmail.com}

\keywords{oscillation; perfect spline; sharp inequalities for derivatives; modulus of continuity of differentiation operator}

\subjclass[2010]{26D10, 41A17, 41A44}
\begin{abstract}
        In this paper we study analogues of the perfect splines for weighted Sobolev classes of functions defined on the half-line. Maximally oscillating splines play important role in the solution of certain extremal problems. In particular, using these splines, we characterize the modulus of continuity of the differential operator. 
\end{abstract}

\maketitle 

\section{Introduction}
Let $I$ be a finite interval or the positive half-line $\RR_+$. 
Let $X$ be a normed space of real-valued functions defined on $I$ and  $f_\pm\in C(I)$ be positive functions. For $x\in X$ set $$\|x\|_{X,f_-,f_+}:=\left\|\frac{\max\{x(\cdot),0\}}{f_+(\cdot)}+\frac{\max\{-x(\cdot),0\}}{f_-(\cdot)} \right\|_X.$$
%If $f_- = f_+ =: f$, we write $\|x\|_{X,f}$ instead of $\|x\|_{X,f_-,f_+}$.
%
%For brevity we denote $\lim\limits_{t\to+\infty}x(t)$ by $x(\infty)$ and write $\|x\|$ instead of $\|x\|_{L_\infty(I)}$.

For positive functions $f_\pm,g\in C(I)$ and natural $r$ set
\begin{multline*}
W_{f_-,f_+,g}^r(I):=\left\{x\in C(I):x^{(r-1)}\in {\rm AC_{loc}},\,\|x\|_{C(I),f_-,f_+}<\infty,\,\right.
\\ \left.\|x^{(r)}\|_{L_\infty(I),g}\leq 1\right\}.
\end{multline*}
%$$W_{f_-,f_+,g}^r(I):=\left\{x\in L_{f_\pm,g}^r(I)  :\|x^{(r)}\|_{L_\infty(I),g}\leq 1\right\}.$$
If $f_-=f_+$, then we write  $W_{f,g}^r(I)$ instead of $W_{f_-,f_+,g}^r(I)$. In the case, when $f\equiv 1$ we write $W_{\infty,g}^r(I)$ instead of $W_{f,g}^r(I)$.

For  $k=1,\ldots,r-1$, we call the function 
\begin{equation}\label{omega}
\omega(W_{f_-,f_+,g}^r(I), D^k,\delta):=\sup\limits_{x\in W_{f_-,f_+,g}^r(I),\, \|x\|_{C(I),f_,f_+}\leq\delta}\|x^{(k)}\|_{C(I)},\,\delta\geq 0
\end{equation}
a modulus of continuity of $k$-th order differentiation operator on the class $W_{f_-,f_+,g}^r(I)$.

%Everywhere below we assume that the functions $f$ and $g$ are positive continuous non-increasing on $\RR_+$.

The study of the function $\omega$ is closely related to the  sharp Landau--Kolmogorov type inequalities. The first results in this topic were obtained in the  1910s by Landau~\cite{Landau13b}, Hadamard~\cite{Hadamard}, Hardy and Littlewood~\cite{HardyLittlewood}. Since then the topic was intensively studied. For a more detailed overview of the history of the question see~\cite{JIA15} and references therein.

Let $I$ denote a segment or a half-line. Suppose that two functions $\Phi, \varphi\in C(I)$ such that $\Phi(t)>\varphi(t)$ for all $t\in I$ are given.
\begin{definition}\label{def::oscillation}
 We say that a function $x\in C(I)$ has $n\in \NN$ points of oscillation between the functions $\varphi$ and $\Phi$, if $\varphi(t)\leq x(t)\leq \Phi(t)$ for all $t\in I$ and there exist points $s_k\in I$, $k=1,\ldots, n$, $s_1<s_2<\ldots<s_{n}$, such that for $k=1,\ldots, n$
\begin{equation}\label{oscillationDef}
x(s_k) = 
\begin{cases}
\Phi(s_k), & k\text { is odd}, \\
\varphi(s_k), & k\text { is even}.
\end{cases}
\end{equation}
\end{definition}

For the case, when $I=[0,\infty)$ and the weights $f_\pm$ and $g$ are constants, the information about the functional $\omega$ follows from the results of Schoenberg and Cavaretta~\cite{schoenbergCavaretta}. We sketch their approach to build the extremal for~\eqref{omega} functions.

One can prove, that for each $A>0$ there exists $\delta = \delta(A,r,n)>0$ and a perfect spline $G = G_{r,n,A}$ of order $r\in\NN$ with $n$ knots, defined on the segment $[0,A]$ that has $n+r+1$ points of oscillation between the constant functions $ -\delta$ and $\delta$. For fixed $r$ and $n$,  $\delta$ is a continuous increasing function of $A$, $\delta\to 0$ as $A\to+0$ and  
\begin{equation}\label{constantNonConstantGDifference}
\delta\to \infty, \text{ as } A\to\infty.
\end{equation}
Moreover, for fixed $r$ and $A$, $\delta$ is a decreasing sequence (numbered by the parameter $n$) and $\delta\to 0$ as $n\to\infty$.

For each $\delta>0$, these properties allow to define a sequence of numbers $A_n$ such that corresponding perfect spline $G_n = G_{r,n,A_n}$ has $n$ knots and oscillates between $-\delta$ and $\delta$ for all $n\in\NN$. Then  $A_n\to \infty$ as $n\to\infty$ and the sequence $G_n$ generates a limiting perfect spline $G_\delta$ defined on the whole half-line with infinite number of knots and  oscillation points between $ -\delta$ and $\delta$. The spline $G_\delta$ is extremal in problem~\eqref{omega}.

In the case, when $g$ is non-constant, the natural substitution for the polynomial splines are so-called $g$-splines.

Let $I = (a,b)$, where $a\in\RR$ and $b$ denotes either a real number or the positive infinity. Let a positive function $g \in C(I)$ be given.

\begin{definition}
The function $G\in C^{r-1}(I)$ is called a perfect $g$-spline of the order $r$ with $n\in \NN$ knots $a<t_1<\ldots<t_n<b$, if on each of the intervals $(t_i,t_{i+1}), i=0,1,\ldots,n$, $t_0:=a, t_{n+1}:=b$, there exists derivative $G^{(r)}$ and $\frac{G^{(r)}(t)}{g(t)}\equiv\epsilon\cdot(-1)^i$ on the intervals $(t_i,t_{i+1}), i=0,1,\ldots,n$, where $\epsilon\in\{1,-1\}$.
\end{definition}

One can repeat the steps described above in the general case of rather arbitrary (we will specify the restrictions for them below) functions $g$ and $f_\pm$, substituting the polynomial perfect splines that oscillate between constant functions, by perfect $g$-splines that oscillate between the functions $-\delta_A\cdot f_-$ and $\delta_A\cdot f_+$, see~\cite{JIA15}, where the symmetric case, when $f_- = f_+$, was considered.

An important difference between the cases of constant and non-constant functions $f_\pm$ and $g$ is in the fact, that~\eqref{constantNonConstantGDifference} may not hold in the latter case. More precisely,  in the case $f_-=f_+$, the following result was proved in~\cite{JIA15}.

Set  $g_0:=g$ and $g_k(t):=\int\limits_0^tg_{k-1}(s)ds$, $t\geq 0$, $k=1,2,\ldots,r$.
\begin{theorem}\label{th::JIA}
	Let $r\in\NN$ and $f,g\in C[0,\infty)$ be non-increasing positive functions. Relation~\eqref{constantNonConstantGDifference} does not hold if and only if
\begin{equation}\label{th0}
	A_0:=\int\limits_0^\infty g(t)dt<\infty,
\end{equation}
for $k=1,\ldots,r-1$
\begin{equation} \label {th1}
		A_k:=\int\limits_0^\infty\left[\sum\limits_{s=0}^{k-1}\frac{(-1)^{k-s-1}A_s}{(k-s-1)!}t^{k-s-1}+(-1)^kg_k(t)\right]dt<\infty,
\end{equation}
and 
\begin{equation}\label{th2}
	\sup\limits_{t\in[0,\infty)}\frac{\left|P_r(t)\right|}{f(t)}<\infty,
\end{equation}
where 
\begin{equation}\label{Pr}
P_r(t) :=  (-1)^r\sum\limits_{s=0}^{r-1}\frac{(-1)^{r-s-1}A_s}{(r-s-1)!}t^{r-s-1}+g_r(t).
\end{equation}
\end{theorem}

In~\cite{JIA15}, the case $f_-=f_+$ was considered, and the following results were obtained.

In the case, when~\eqref{constantNonConstantGDifference} holds, using Schoenberg and Cavaretta's  approach, the values $\omega(W_{f,g}^r[0,\infty), D^k,\delta)$ were characterized in terms of the limit $g$-splines with infinite number of knots and points of oscillation for all $\delta>0$.

In the case, when~\eqref{constantNonConstantGDifference} does not hold, the values of the functional $\omega(W_{f,g}^r[0,\infty), D^{k},\delta)$ were characterized only for a non-increasing sequence of numbers $\left\{\delta_{r,n}\right\}_{n=0}^{\infty}$. Moreover, this sequence may not tend to $0$ as $n\to\infty$.

In this article we characterize the values of $\omega(W_{f_-,f_+,g}^r[0,\infty), D^{k},\delta)$ for all $\delta>0$ in the case when~\eqref{constantNonConstantGDifference} does not hold. However, we impose stricter restrictions on the functions $f_\pm$ and $g$, then those described in Theorem~\ref{th::JIA}. Namely, we characterize $\omega(W_{f_-,f_+,g}^r[0,\infty), D^{k},\delta)$ under the following three assumptions.

\begin{assumption}\label{a::g}
The function $g\in C[0,\infty)$ is positive non-increasing and such that conditions~\eqref{th0} and~\eqref{th1} hold.
\end{assumption}

\begin{assumption}\label{a::fpm}
The functions $f_\pm\in C[0,\infty)$ are non-increasing positive and such that 
\begin{equation}\label{f>0}
f_\pm(\infty) >0
\end{equation}
\end{assumption}

\begin{assumption}\label{a::fg}
\begin{equation}\label{assumption3}
\varliminf\limits_{t\to\infty}\frac {f_\pm(t)-f_\pm(\infty)}{|P_r(t)|} =0,
\end{equation}
where the function $P_r$ is defined in~\eqref{Pr}.
\end{assumption}

Note, that Assumption~\ref{a::g} is the same as in Theorem~\ref{th::JIA}, but the assumption on the functions $f_\pm$ are stricter than the ones in Theorem~\ref{th::JIA}; in particular Assumptions~\ref{a::g} and~\ref{a::fpm} imply~\eqref{th2} in the case $f_-=f_+$. However, the important case, when $f_\pm$ are constant, is included.

According to Definition~\ref{def::oscillation}, the points of oscillation are 'positively orientated', that is $x(s_1) = \Phi(s_1)$. We can define 'negatively orientated' points of oscillation, substituting  condition~\eqref{oscillationDef} by 
$$
x(s_k) = 
\begin{cases}
\Phi(s_k), & k\text { is even}, \\
\varphi(s_k), & k\text { is odd}.
\end{cases}
$$
\begin{remark}\label{r::positiveVsNegativeOcsillationOrientation}
Let a function $x_+ = x_+(\varphi,\Phi)$ have $n$ positively orientated points of oscillation between the functions $\varphi$ and $\Phi$. Then the function $x_- =-x_+(-\Phi,-\varphi)$ has $n$ negatively orientated points of oscillation between the functions $\varphi$ and $\Phi$.

In particular, if $\Phi$ is positive and $\varphi\equiv -\Phi$, then $x_- = -x_+$. However, there is no such relation in general situation.   
\end{remark}

\begin{definition}
	We call a function $x\in C[0,\infty)$ $n$-piecewise monotone, $n\in\NN$, if there exists $\epsilon\in \{1, -1\}$ and $0=:t_0<t_1<\ldots < t_{n-1} < t_n:=\infty$, so that $\epsilon\cdot (-1)^k x$ is increasing on the interval $(t_{k-1}, t_{k})$, $k=1,\ldots, n$.
\end{definition}

For a locally absolutely continuous function $x$ defined on an interval, the notations  ${\rm sgn\,}x'(t) = 1$ (${\rm sgn\,}x'(t) = -1$) will mean that the function $x$ is increasing (decreasing) in some neighborhood of the point $t$.

\begin{definition}
A primitive $G$ of the order $r$ of the function $g$ or $-g$ on $I$ will be called a perfect $g$-spline of the order $r$ with $0$ knots.
\end{definition}

\begin{definition}
Denote by $\Gamma^{r}_{n,g}(I)$ the set of all perfect $g$-splines $G$ defined on $I$ of the order $r$ with not more than $n\in \ZZ_+$ knots.
\end{definition}

To obtain our main result, we study some properties of the maximally oscillating perfect $g$-splines. 
We prove the following result, which also has an independent interest. In particular, statement~(c) of the theorem states, that maximally oscillating splines are least deviating from zero in non-symmetric weighted norm among the $g$-splines of the class $\Gamma^{r}_{n,g}(I)$. Such property for polynomial splines is well known.
\begin{theorem}\label{th::ZolotarevSplines}
Let Assumptions~\ref{a::g}, \ref{a::fpm}, and~\ref{a::fg}  hold. 
 \begin{enumerate}
 \item[(a)] There exist two non-increasing sequences $\left\{ a_n^\pm\right\}_{n=1}^\infty$ of positive numbers  such that $\lim\limits_{n\to\infty}a_n^\pm = 0$, and for each $\tau\in [a_{n+1}^\pm, a_n^\pm)$, $n\in\ZZ_+$, $a_0^\pm:=\infty$, there exists a perfect $g$-spline $G^\pm_\tau\in \Gamma^r_{n,g}[0,\infty)$ with exactly $n$ knots and $n+1$ (positively oriented for $G_\tau^+$ and negatively oriented for $G_\tau^-$) points of oscillation between $-\tau f_-$ and $\tau f_+$.
 \item[(b)] For all $\tau \in [a^\pm_{n+1}, a^\pm_n)$ and $s=0,1,\ldots, r-1$, $(G^\pm_\tau)^{(s)}$ is $(n+1)$-piecewise monotone. Moreover, for $s = 0,\ldots, r$,
 \begin{equation}\label{gSplineDerivativesAtZero}
 	{\rm sgn\,}(G^\pm_\tau)^{(s)}(0) = \pm (-1)^s.
 \end{equation}

\item[(c)] For each $n\in\ZZ_+$
$$\inf\limits_{G\in \Gamma^r_{n,g}[0,\infty)}\|G\|_{C[0,\infty), f_-,f_+} = \min\{a_{n+1}^+, a_{n+1}^-\}.
$$
If $a_{n+1}^+< a_{n+1}^-$, then the infimum is attained on $G_{a_{n+1}^+}^+$, otherwise it is attained on $G_{a_{n+1}^-}^-$. 
\end{enumerate}
\end{theorem}

\begin{remark}
It will follow from the proof, that the $g$-splines $G_{a_{n+1}^\pm}^\pm$ satisfy
\begin{equation}\label{g+AtInfty}
G_{a_{n+1}^+}^+(\infty)
=
\begin{cases}
 -a_{n+1}^+f_-(\infty), & n\,\text{is even,} \\
a_{n+1}^+f_+(\infty), & n\,\text{is odd.}
\end{cases}
\end{equation}
and 
\begin{equation}\label{g-AtInfty}
G_{a_{n+1}^-}^-(\infty)
=
\begin{cases}
 -a_{n+1}^-f_-(\infty), & n\,\text{is odd,} \\
a_{n+1}^-f_+(\infty), & n\,\text{is even.}
\end{cases}
\end{equation}
These $g$-splines can be viewed as the analogues of the perfect Chebyshev splines, having an 'oscillation point at $\infty$'. Here and below $f(\infty)$ denotes $\lim\limits_{t\to+\infty}f(t)$.
\end{remark}

In the theorem above, the perfect $g$-splines are indexed by their non-symmetric weighted norms. The same splines admit indexation based on their values at infinity. We formulate the result for the splines with positively oriented oscillation points; a similar results holds in the case of negatively oriented oscillation points.
\begin{theorem}\label{th::InfinityValues}
Let Assumptions~\ref{a::g} and \ref{a::fpm} hold and $n\in\ZZ_+$ be given. For all 
$$
\alpha\in 
\begin{cases}
[0,\infty),& n\text{ is odd}, \\
(0,\infty],& n\text{ is even},
\end{cases}
$$
there exists $\tau=\tau_n(\alpha)>0$
and a perfect $g$-spline $G=G(n,\alpha)\in \Gamma^r_{n,g}[0,\infty)$ with exactly $n$ knots that has $n+1$ positively oriented points of oscillation between $-\tau\cdot f_-$ and $\tau\cdot f_+$ and such that  
\begin{equation}\label{infiniteBoundaryCondition1}
\frac{\tau f_+(\infty) - G(\infty)}{\tau f_-(\infty) + G(\infty)} = \alpha,
\end{equation}
where~\eqref{infiniteBoundaryCondition1} is understood as $ G(\infty) = - \tau f_-(\infty)$ if $\alpha=\infty$.
For $n\in\ZZ_+$, the function $\tau_n$ is continuous, non-decreasing in the case of odd $n$, and is non-increasing in the case of even $n$.
For all $n\in \NN$ and $\alpha,\beta >0$ one has $\tau_n(\alpha)\leq \tau_{n-1}(\beta)$.
\end{theorem}
\begin{remark}
We do not prove the uniqueness of the maximally oscillating perfect $g$-splines in Theorem~\ref{th::ZolotarevSplines} and~\ref{th::InfinityValues}. However, the correctness of the definition of the functions $\tau_n$ in Theorem~\ref{th::InfinityValues} will be proved below.
\end{remark}
\begin{remark} Note, that Assumption~\ref{a::fg} is needed in Theorem~\ref{th::ZolotarevSplines}, but not needed in Theorem~\ref{th::InfinityValues}. It guarantees that there is no 'gap' between the possible values of the non-symmetric weighted norms of the perfect $g$-splines with $n-1$ and $n$ knots, $n\in\NN$.
\end{remark}

The main result of this article is the following theorem.
\begin{theorem}\label{th::modulusOfContinuity}
Let $r\in\NN$, $r\geq 2$ and Assumptions~\ref{a::g}, \ref{a::fpm}, and~\ref{a::fg}  hold. 
Then for all $k =1,\ldots, r-1$ and $\delta >0$, 
$$\omega(W_{f_-,f_+,g}^r[0,\infty), D^k,\delta) = \max\left\{|(G_\delta^+)^{(k)}(0)|,|(G_\delta^-)^{(k)}(0)|\right\},$$
where $G_\delta^\pm$ are perfect $g$-splines from Theorem~\ref{th::ZolotarevSplines}. If 
$$|(G_\delta^+)^{(k)}(0)| > |(G_\delta^-)^{(k)}(0)|,$$ then the supremum in the definition of $\omega(W_{f_-,f_+,g}^r[0,\infty), D^k,\delta)$ is attained on $G_\delta^+$, otherwise it is attained on $G_\delta^-$.
\end{theorem}
\begin{remark}
Note, that if $\delta \geq \max\{a_1^+, a_1^-\}$, then $G_\delta^\pm$ have $0$ knots, and hence $G_\delta^\pm = \pm|P_r| + C_\delta^\pm$, $C_\delta^\pm\in \RR$. This implies that  the functional $\omega(W_{f_-,f_+,g}^r[0,\infty), D^k,\cdot)$ becomes constant for all big enough values of $\delta$.
\end{remark}

The article is organized as follows. In Section~2 we  prove the existence of maximally oscillating $g$-splines on the finite segments and the half-line. In Section~3 we study the non-symmetrical weighted norms of the maximally oscillating $g$-splines. Section~4 is devoted to the proofs of the main results.

%\begin{theorem}[Brouwer~\cite{brouwer}]
%For any continuous function $\phi\colon K\to K$ mapping a compact convex set $K\subset \RR^n$ to itself, there is a point $\xi^*$ such that $\phi(\xi^*) = \xi^*$.
%\end{theorem}

\section{Perfect $g$-splines}
\subsection{Auxiliary definitions}
In this article, we will often count or estimate the number of zeros and sign changes of the functions. Below we adduce necessary definitions.

Let $x$ be a continuous on an interval function and $ t_1<\ldots<t_n$ be its zeros, $n\geq 2$. The points $t_1,\ldots, t_n$ are called separated zeros, if the function $x$ is not identical zero on each of the intervals $(t_i,t_{i+1})$, $i=1,\dots, n-1$. 

We say that a function $x$ has essential sign change on the interval $I$, if both sets $\{t\in I\colon x(t) >0\}$ and $\{t\in I\colon x(t) <0\}$ have positive measures. 

We say that a function $x$ has essential change of sing at the point $z$, if there exists $\varepsilon>0$ such that for almost all $0<u<\varepsilon$, ${\rm sgn\,}x(z-u) = -  {\rm sgn\,}x(z+u)\neq 0$. 

We say that a function $x$ has exactly $n$ essential sign changes on the interval $I$, if there exist $n+1$ points $t_1<\ldots <t_{n+1}$ from this interval and a number $\varepsilon>0$ such that for almost all $u\in (-\varepsilon,\varepsilon)$ one has
${\rm sgn\,} x(t_i+u)= -{\rm sgn\,} x(t_{i+1}+u)\neq 0$, $i=1,\ldots, n$, and there exist no system of $n+2$ points with such property.

For continuous functions we use the term 'sign change' instead of 'essential sign change'.
%If $x$ is continuous on the interval $I$ function and has $n$ sign changes, then $x$ has at least $n$ separated zeros on $I$.

We will often use the following analogue of the Rolle theorem.

Between two separated zeros of a locally absolutely continuous function $x$, there is an essential change of sign of the function $x'$.

\subsection{Perfect $g$-splines with given zeros}
We need the following topological result. 
\begin{theorem}[Borsuk~\cite{borsuk}]
Let $S^n = \{\xi\in \RR^{n+1}\colon \|x\| = r\}$, where $\|\cdot\|$ is some norm in $\RR^{n+1}$, and $\phi\colon S^{n}\to\RR^n$ be a continuous odd function. Then there exists a point $\xi^*\in S^n$ such that $\phi(\xi^*) = 0$.
\end{theorem}
According to Velikin~\cite{velikin}, the idea to use the Borsuk theorem in the proof of existence of perfect splines with given zeros, belongs to Ruban.
\begin{lemma}\label{l::splineWithGivenZerosExistance}
Let $A>0$, $n\in\NN$ and $0\leq s_1<s_2<\ldots<s_n<A$. There exists a perfect $g$-spline $G\in \Gamma^{r}_{n,g}[0,A]$ such that 
\begin{equation}\label{boundaryConditions}
G^{(k)}(A) = 0,\, k=0,1,\ldots, r-1.
\end{equation}
and 
\begin{equation}\label{givenZeros}
G(s_k) = 0, k=1,\ldots, n.
\end{equation}
Moreover, conditions~\eqref{boundaryConditions} and~\eqref{givenZeros} imply that $G$ changes sign in each of the points $s_k$, $k=1,\ldots, n$, and has exactly $n$ knots.
\end{lemma}

Consider the sphere $S^{n}:=\left\{(\xi_1,\ldots, \xi_{n+1})\in\RR^{n+1}\colon \sum\limits_{k=1}^{n+1}|\xi_k| = A\right\}$. For arbitrary $\xi\in S^n$ consider the partition of the segment $[0,A]$ by the points $t_m = \sum\limits_{k=1}^m|\xi_k|$, $m=1,\ldots, n$. Let $G(\xi)\in C^{r-1}[0,A]$ be the function that satisfies boundary conditions~\eqref{boundaryConditions} and such that $G^{(r)}(\xi,t) = {\rm sgn\,}\xi_k\cdot g(t)$ on the interval $(t_{k-1},t_k)$, $k=1,\ldots, n+1$, $t_0:=0$, $t_{n+1} = A$; such function is uniquely determined by the imposed conditions. Moreover, $G(\xi) = - G(-\xi)$ for arbitrary $\xi\in S^n$ and $G(\xi)$ uniformly converges to $G(\xi_0)$ provided by $\xi\to\xi_0$.

Consider a map $\phi\colon S^n\to\RR^n$, $\phi(\xi) = (G(\xi;s_1),\ldots, G(\xi;s_n))$. It is continuous and odd. By the Borsuk theorem, there exists $\xi^*\in S^n$ such that $\phi(\xi^*) = 0$. The function $G(\xi^*)$ satisfies~\eqref{boundaryConditions} and~\eqref{givenZeros},  $G^{(r)}(\xi^*)$ is non-zero almost everywhere. Hence all zeros $s_1,\ldots, s_n$ of $G(\xi^*)$ are separated, and thus the function $G^{(r)}(\xi^*)$ has at least $n$ sign changes due to the Rolle theorem. From the construction it follows, that the function $G^{(r)}(\xi^*)$ can not have more than $n$ sign changes, hence $G(\xi^*)\in \Gamma^{r}_{n,g}[0,A]\setminus \Gamma^{r}_{n-1,g}[0,A]$. All zeros $s_1,\ldots, s_n$ are simple, since otherwise the Rolle theorem would imply existence of more than $n$ sign changes of the function $G^{(r)}(\xi^*)$. Hence $G(\xi^*)$ is a desired perfect $g$-spline. The lemma is proved.

\begin{lemma}\label{l::splineWithGivenZerosUniqueness}
Let $A>0$, $n\in\NN$ and $0 < s_1<s_2<\ldots<s_n<A$. If two perfect $g$-splines $G_1, G_2\in \Gamma^{r}_{n,g}[0,A]$ satisfy conditions~\eqref{boundaryConditions} and~\eqref{givenZeros}, then either $G_1 \equiv G_2$, or $G_1 \equiv - G_2$.
\end{lemma}

We can assume that both $G_1$ and $G_2$ are positive on $[0,s_1)$. Hence $G_1$ and $G_2$ have same signs on each of the intervals $(s_k,s_{k+1})$, $k=0,\ldots, n$, $s_0:= 0$, $s_{n+1}:=A$. 

Let $k\in \{0,\ldots, n\}$ and $s\in (s_k,s_{k+1})$ be fixed. We prove that 
\begin{equation}\label{equalityInS}
G_1(s) = G_2(s). 
\end{equation}
Assume the contrary, without loss of generality we may assume that $|G_1(s)| < |G_2(s)|$. There exists $\varepsilon\in (0,1)$ such that $G_1(s) = (1-\varepsilon)\cdot G_2(s)$. The difference $G:=G_1 - (1-\varepsilon)\cdot G_2$ satisfies conditions~\eqref{boundaryConditions} and~\eqref{givenZeros}; moreover, it has an additional zero at the point $s$. From the definition of the function $G$ it follows that the function $G^{(r)}$ is non-zero almost everywhere and changes its sign only at the knots of the perfect $g$-spline $G_1$. Hence all zeros $s_1,\ldots, s_n$, $s$ and $A$ of the function $G$ are separated and thus the Rolle theorem implies that the function $G^{(r)}$ has at least $n+1$ sign changes. This contradiction proves~\eqref{equalityInS}. Due to arbitrariness of the point $s$, this implies that $G_1\equiv G_2$.
The lemma is proved.

\begin{lemma}\label{l::continuityOfKnotsZerosDependence}
Let $A>0$ and $n\in\NN$ be given. Suppose that for each $m=0,1,\ldots$ a point  $ s^{(m)} = (s_{1,m}, \ldots, s_{n,m})\in\RR^n$,  $0\leq s_{1,m}<\ldots<s_{n,m}<A$ is given. Let $G_m\in \Gamma^{r}_{n,g}[0,A]$ be a perfect $g$-spline that satisfies~\eqref{boundaryConditions} and vanishes at the points $s^{(m)}$, $m=0,1,\ldots$. Denote by $t^{(m)} = (t_{1,m},\ldots, t_{n,m})$, $0< t_{1,m}<\ldots<t_{n,m}<A$ the knots of $G_m$. Then $s^{(m)}\to s^{(0)}$ as $m\to\infty$ implies $t^{(m)}\to t^{(0)}$ as $m\to\infty$.
\end{lemma}

Assume, that the sequence $\left\{t^{(m)}\right\}_{m=1}^\infty$ has two different limit points $u=(u_1,\ldots, u_n)$ and $v = (v_1,\ldots, v_n)$. Consider the perfect $g$ splines $G_u$ and $G_v$ with knots at the point $u$ and $v$ that satisfy boundary conditions~\eqref{boundaryConditions}; each of the splines is determined up to the sign. Moreover, since the perfect spline continuously  (in the sense of uniform convergence) depends on its knots, both $G_u$ and $G_v$ vanish at the points $s^{(0)}$. However, this contradicts to Lemma~\ref{l::splineWithGivenZerosUniqueness}. Hence the sequence  $\left\{t^{(m)}\right\}_{m=1}^\infty$ has a limit. It can't be different from $t^{(0)}$ due to Lemma~\ref{l::splineWithGivenZerosUniqueness}. The lemma is proved.

\subsection{Maximally oscillating perfect $g$-splines on a segment}
\begin{remark}
Everywhere below for brevity we write 'a function has $n$ points of oscillation' instead of 'a function has $n$ positively orientated points of oscillation'.
\end{remark}

The technics that involve Brouwer fixed-point theorem in the proof of existence of oscillating functions, can be found for example in~\cite[\S 10, Chapter 2]{karlinStudden}.

\begin{lemma}\label{l::finiteIntervalOscillation}

Let $A>0$ and two functions $f_+,f_-\in C[0,A]$ be given. Assume $f_\pm(t) > 0$ for all $t\in [0,A]$. Then for each $n\in \ZZ_+$ there exists $\tau>0$ and a perfect $g$-spline $G\in  \Gamma^{r}_{n,g}[0,A]$ that satisfies~\eqref{boundaryConditions}, has exactly $n$ knots, and has $n+1$ points of oscillation between the functions $-\tau\cdot f_-$ and $\tau\cdot f_+$.
\end{lemma}

For all $\varepsilon\in \left(0,\frac A{n+1}\right)$ consider the $n$-dimensional simplex 
$$
\Xi^n_\varepsilon :=\left\{(\xi_1,\ldots, \xi_{n+1})\in\RR^{n+1}\colon \xi_1,\ldots, \xi_{n+1}\geq\varepsilon, \sum\limits_{k=1}^{n+1}\xi_k = A\right\}.$$
Set
$$
\Xi^n_0 :=\left\{(\xi_1,\ldots, \xi_{n+1})\in\RR^{n+1}\colon \xi_1,\ldots, \xi_{n+1}>0, \sum\limits_{k=1}^{n+1}\xi_k = A\right\}.$$

For each point $\xi\in \Xi^n_\varepsilon$ consider the partition of the segment $[0,A]$ by the points $s_m = \sum\limits_{k=1}^m|\xi_k|$, $m=1,\ldots, n$. Due to the definition of the simplex, $0<s_1<\ldots< s_n<A$. Hence, in virtue of Lemma~\ref{l::splineWithGivenZerosExistance}, there exists a perfect $g$-spline $G(\xi)$ with exactly $n$ knots 
%$0<t_1(\xi)<\ldots<t_n(\xi)<A$ 
that satisfies~\eqref{boundaryConditions} and~\eqref{givenZeros}. The sign of $G(\xi)$ is chosen in such a way, that $G(\xi)$ is positive on $(0,s_1)$. For each $k = 1,\ldots, n+1$, set 
$$
\delta_k(\xi) := \min\left\{\delta \geq 0\colon |G(\xi;t)|\leq \delta\cdot f_k(t),t\in (s_{k-1},s_{k})\right\},
$$
where $s_0:=0$, $s_{n+1}:=A$, $f_k = f_+$ for odd $k$, and  $f_k = f_-$ for even $k$.

Set $\Delta_k(\xi):= \delta_k(\xi) - \min\limits_{j=1,\ldots,n+1}\delta_j(\xi)$, $k=1,\ldots,n+1$, and $\Delta(\xi):=\sum\limits_{k=1}^{n+1}\Delta_k(\xi)$. 

If for some $\varepsilon>0$ there exists $\xi\in\Xi_\varepsilon^{n}$ such that $\Delta(\xi)=0$, then corresponding perfect $g$-spline $G(\xi)$ is a desired one.

Assume that for all $\varepsilon>0$ and $\xi\in \Xi_\varepsilon^n$ 
\begin{equation}\label{noOscillationAssumption}
\Delta(\xi)\neq 0.
\end{equation}

In virtue of Lemma~\ref{l::continuityOfKnotsZerosDependence}, $\delta_k$ and $\Delta_k$ are continuous functions of $\xi$, $k=1,\ldots,n+1$, hence $\Delta$ is also a continuous function of $\xi$.

Note, that since the functions $f_+$ and $f_-$ are separated from zero, 
$$\sup\limits_G\|G'\|_{C[0,A]}<\infty,$$
where the supremum is taken over all splines ${G\in\Gamma^{r}_{n,g}[0,A]}$ that satisfy condition~\eqref{boundaryConditions}. From this fact  it follows, that  there exists $\alpha>0$ that does not depend on $\xi$ and $\varepsilon$ such that 
\begin{equation}\label{delta<xi}
\delta_k(\xi)\leq\alpha\cdot\xi_k,\,k=1,\ldots, n+1.
\end{equation}
From $\inf\limits_{G\in\Gamma^{r}_{n,g}[0,A]}\|G\|_{C[0,A]} > 0$ it follows that
\begin{equation}\label{maxDeltaIsSeparatedFromZero}
\inf\limits_{\varepsilon >0}\min\limits_{\xi\in\Xi_\varepsilon^n}\max\limits_{k=1,\ldots, n+1}\delta_k(\xi) =:\beta >0.
\end{equation}
We claim that
\begin{equation}\label{deltaIsSeparatedFromZero}
\inf\limits_{\varepsilon >0}\min\limits_{\xi\in\Xi_\varepsilon^n}\Delta(\xi)=:\gamma >0.
\end{equation}
Indeed, if~\eqref{deltaIsSeparatedFromZero} does not hold, then, due to~\eqref{noOscillationAssumption}, for arbitrary $\varepsilon>0$ there exists $\xi_\varepsilon\in \Xi_0^n\setminus \Xi_\varepsilon^n$ such that $\Delta(\xi_\varepsilon) < \frac \beta 2$. Hence 
$$\frac \beta 2>\Delta(\xi_\varepsilon)\geq \max_{k=1,\ldots n+1}\Delta_k(\xi_\varepsilon) = \max_{k=1,\ldots n+1}\delta_k(\xi_\varepsilon) - \min\limits_{k=1,\ldots, n+1}\delta_k(\xi_\varepsilon),$$ and due to~\eqref{maxDeltaIsSeparatedFromZero} this implies that $\min\limits_{k=1,\ldots, n+1}\delta_k(\xi_\varepsilon) > \frac \beta 2$. However, this contradicts to~\eqref{delta<xi}.

Define the map $\phi\colon\Xi_\varepsilon^n\to\Xi_\varepsilon^n$, 
$$
\phi(\xi) =  \varepsilon + \frac{A-(n+1)\varepsilon}{\Delta(\xi)}\left(\Delta_2(\xi),\Delta_3(\xi),\ldots,\Delta_{n+2}(\xi)\right),
$$
where $\Delta_{n+2}(\xi):=\Delta_{1}(\xi)$. By the construction, $\phi$ is continuous, hence by the Brouwer fixed point theorem there exists $\xi^*\in\Xi_\varepsilon^n$ such that 
\begin{equation}\label{fixedPoint}
\phi(\xi^*) = \xi^* = (\xi^*_1,\ldots, \xi^*_{n+1}).
\end{equation}
Hence for $k=1,\ldots, n+1$, due to~\eqref{delta<xi} and~\eqref{deltaIsSeparatedFromZero}, one has 
\begin{equation}\label{estimateForXi}
\xi^*_k =\varepsilon + \frac{A-(n+1)\varepsilon}{\Delta(\xi^*)}\Delta_{k+1}(\xi^*) < \varepsilon + \frac A \gamma\delta_{k+1}(\xi^*) < \varepsilon + \frac {\alpha\cdot A} \gamma \xi^*_{k+1}.
\end{equation}
From the construction of the functions $\Delta_k$ it follows that there exists $k\in \{1,\ldots, n+1\}$ such that $\Delta_k(\xi^*) = 0$. Let for definiteness $\Delta_1(\xi^*) = 0$. Then, due to~\eqref{fixedPoint}, $\xi^*_{n+1} = \varepsilon$. Consecutive application of estimate~\eqref{estimateForXi} implies that there exist positive numbers $\eta_1,\ldots, \eta_n$ that are independent of $\varepsilon$ and $\xi^*$ and such that $\xi^*_n<\eta_n\varepsilon$, $\xi^*_{n-1}<\eta_{n-1}\varepsilon$ and so on, $\xi^*_{1}<\eta_{1}\varepsilon$. However, if $\varepsilon$ is chosen to be enough small, this contradicts to the fact that $\sum\limits_{k=1}^{n+1}\xi^*_k = A$. The lemma is proved.

\subsection{Auxiliary results}
We write $W^r_{\infty,\infty}[0,\infty)$ instead of $ W^r_{\infty,g}[0,\infty)$, if $g\equiv 1$.
The following result is essentially due to Landau~\cite{Landau13a}. We will give its proof for completeness.
\begin{lemma}\label{l::derivativeAtInfy}
Let $r\in \NN$, $r\geq 2$ and $x \in W^r_{\infty,\infty}[0,\infty)$ be given. If $\lim\limits_{t\to\infty}x(t)$ exists, then $\lim\limits_{t\to\infty}x^{(k)}(t) = 0$, $k=1,\ldots, r-1$.
\end{lemma}

The function $x^{(k)}$ is bounded, since $x$ and $x^{(r)}$ are bounded, $k=1,\ldots, r-1$. The lemma can be proved by induction on $k$.

Let $M>0$, $t>M$ and $\varepsilon >0$. Then
$x(t+\varepsilon) - x(t) =\varepsilon x'(t) + \frac{\varepsilon^2} 2 x''(\xi),$
$\xi\in (t,t+\varepsilon)$, hence 
\begin{equation}\label{x'}
|x'(t)|\leq  \frac {|x(t+\varepsilon) - x(t)|}{\varepsilon} + \frac \varepsilon 2 |x''(\xi)|.
\end{equation}
 Since $\lim\limits_{t\to\infty}x(t)$ exists, there exists $M=M(\varepsilon)>0$ such that $|x(s) - x(t)|<\varepsilon^2$ for all $t,s>M$. Hence, due to the fact, that $x''$ is bounded, the right-hand side of~\eqref{x'} can be made arbitrarily small. This implies the statement of the lemma in the case $k=1$. 
 
 Induction step can be done using the same arguments. The lemma is proved. 

\begin{lemma}\label{l::derivativeAtInfy2}
Let $r\in \NN$, $r\geq 2$ and $x \in W^r_{\infty,\infty}[0,\infty)$ be given. If $x^{(r)}$ has a finite number of essential sign changes, then the limit $\lim\limits_{t\to\infty}x(t)$ exists and  $\lim\limits_{t\to\infty}x^{(k)}(t) = 0$, $k=1,\ldots, r-1$.
\end{lemma}

Due to conditions of the lemma, $x'$ has a finite number of sign changes. Hence $x$ is monotone in some neighborhood of infinity. Since it is bounded, the limit $x(\infty)$ exists. To finish the  proof, it is enough to apply Lemma~\ref{l::derivativeAtInfy}. The lemma is proved.

\begin{lemma}\label{l::singAlternation}
Let $x \in W^r_{\infty,\infty}[0,\infty)$ be such, that $x^{(r)}$ is non-zero almost everywhere and has at most $n\in \ZZ_+$ essential sign changes. Assume $x'$ has at least $n$ sign changes. Then each of the functions $x^{(k)}$, $k=0,\ldots, r-1$, is $(n+1)$-piecewise monotone. Moreover, if $\epsilon\in\{1,-1\}$ and $\epsilon\cdot x$ is increasing on the first monotonicity interval, then 
\begin{equation}\label{signsAtZero}
{\rm sgn\,} x^{(k)}(0) = \epsilon (-1)^{k+1}, k = 1,\ldots, r.
\end{equation}

\end{lemma} 

From the conditions of the lemma it follows, that the function $x$ and all its derivatives are non-zero almost everywhere. By Lemma~\ref{l::derivativeAtInfy2}, $x'(\infty) = 0$. Hence $x''$ has at least $n$ sign changes and, by Lemma~\ref{l::derivativeAtInfy2}, $x''(\infty) = 0$. Continuing the same way, we obtain that $x^{(r)}$ has at least $n$ essential sign changes. Since $x^{(r)}$ has at most essential $n$ sign changes, each of the derivatives $x^{(k)}$ has exactly $n$ sign changes, $k=1,\ldots, r$.  Then each of the functions $x^{(k)}$, $k=0,\ldots, r-1$, is $(n+1)$-piecewise monotone. Moreover, the function $x^{(k)}$, $k=1,\ldots, r-1$, changes sign on each of its monotonicity interval, except the last one (where it tends to $0$). This implies~\eqref{signsAtZero}. The lemma is proved.

\begin{definition}\label{def::Pk}
Let Assumption~\ref{a::g} hold. Set $P_0:=g$ and for $k=1,\ldots, r$ and $t\in [0,\infty)$ set 
$$P_k(t) :=  (-1)^k\sum\limits_{s=0}^{k-1}\frac{(-1)^{k-s-1}A_s}{(k-s-1)!}t^{k-s-1}+g_k(t).$$
Note, that $P_k$ is the $k$-th order primitive of the function $g$ such that $P_k(\infty) = 0$, $k =1,\ldots, r$. Moreover, $P_k$ is monotone and does not change sign on $[0,\infty)$.
\end{definition}

\begin{lemma}\label{l::GLessP}
Let Assumption~\ref{a::g} hold, $r\in \NN$, $x\in W^r_{\infty,g}[0,\infty)$ and $\lim\limits_{t\to\infty}x(t) = A\in \RR$. Then for all $t\in \RR$, $|x(t)-A|\leq |P_r(t)|$.
\end{lemma}

It is enough to prove the lemma in the case when $A = 0$. This, in turn, can be done by induction on $r$. Really, using induction hypothesis for $r>1$ (which can be applied, due to Lemma~\ref{l::derivativeAtInfy}), or the definition of the class $W_{\infty, g}^r[0,\infty)$ for $r=1$, we have
\begin{gather*}
|x(t)| = |x(t) - x(\infty)| = \left|\int\limits_{t}^{\infty}x'(s)ds\right| \leq \int\limits_{t}^{\infty}\left|x'(s)\right|ds
\leq \int\limits_{t}^{\infty}\left|P_{r-1}(s)\right|ds 
\\ =
\left|\int\limits_{t}^{\infty} P_{r-1}(s)ds\right| 
=|P_r(t)- P_r(\infty)| = |P_r(t)|
.
\end{gather*}
The lemma is proved.
\subsection{Maximally oscillating perfect $g$-splines on a half-line}

In what follows we assume that Assumptions~\ref{a::g} and~\ref{a::fpm}  hold.

\begin{lemma}\label{l::halfLineOscillation}
Let $n\in\ZZ_+$ be given. For all $\alpha\in (0,\infty)$ there exists $\tau>0$
and a perfect $g$-spline $G=G(n,\alpha)\in \Gamma^r_{n,g}[0,\infty)$ with exactly $n$ knots that has $n+1$ points of oscillation between $-\tau\cdot f_-$ and $\tau\cdot f_+$ and such that  
\begin{equation}\label{infiniteBoundaryCondition}
\frac{\tau f_+(\infty) - G(\infty)}{\tau f_-(\infty) + G(\infty)} = \alpha.
\end{equation}

For all $0<\alpha_1 < \alpha_2<\infty$ there exists a number $C = C(\alpha_1, \alpha_2)>0$ such that   
$\tau\leq C$, 
provided by $\alpha\in [\alpha_1,\alpha_2]$.
\end{lemma}

Set $a:= \frac {1}{1+\alpha} f_+(\infty) - \frac {\alpha}{1+\alpha} f_-(\infty) $. Then $a\in (-f_-(\infty), f_+(\infty))$ and hence both $f_+ -a$ and $f_-+a$ are positive continuous non-increasing on $[0,\infty)$ functions. 

For arbitrary $A>0$ consider a perfect $g$-spline $G_A\in \Gamma_{n,g}^r[0,A]$ that oscillates maximally between the functions $-\tau_A(f_-+a)$ and $\tau_A(f_+ -a)$, where $\tau_A>0$ and the existence of such spline is guaranteed by Lemma~\ref{l::finiteIntervalOscillation}. Let $0\leq s_1^A<\ldots < s_{n+1}^A \leq A$ be the oscillation points of $G_A$ and $0< t_1^A<\ldots < t_{n}^A < A$ be its knots. 

Letting $A\to\infty$ and switching to a subsequence, if needed, we may assume that each of the sequences $s_k^A$, $k=1,\ldots, n+1$, $t_k^A$, $k=1,\ldots, n$ and $\tau_A$, has a finite or infinite limit $s_k$, $t_k$ and $\tau$ respectively.

Since $G_A$ satisfies~\eqref{boundaryConditions}, it can be continued to a function from the class $W^r_{\infty,g}[0,\infty)$ by setting it equal to $0$ on $[A,\infty)$. Then, due to Lemma~\ref{l::GLessP}, for all $A>0$ and $t\in [0,A]$
\begin{equation}\label{G_AIsBounded}
|G_A(t)| \leq |P_r(t)|.
\end{equation}

The latter inequality implies that for all $A>0$ one has $\tau_A \leq \|P_r\|_{C[0,\infty),f_-+a,f_+-a}$. Hence $\tau\leq \|P_r\|_{C[0,\infty),f_-+a,f_+-a}$. If $0<\alpha_1 < \alpha_2<\infty$ are fixed, then both functions $f_-+a$ and $f_+-a$ are separated from zero by a constant independent of $\alpha\in [\alpha_1,\alpha_2]$ due to~\eqref{f>0}. Hence $\|P_r\|_{C[0,\infty),f_-+a,f_+-a}$ is uniformly bounded from above for such values of $\alpha$.

Note, that  $\tau>0$, since otherwise for each $\varepsilon>0$ we could find $A=A(\varepsilon)>1$, such that $\|G_A\|_{C[0,1]}<\varepsilon$, which is impossible, since the restriction of the $g$-spline $G_A$ to $[0,1]$ belongs to $\Gamma^r_{n,g}[0,1]$ and $\inf\limits_{G\in \Gamma^r_{n,g}[0,1]}\|G\|_{C[0,1]}>0$.

Since~\eqref{G_AIsBounded} holds for $t= s_{n+1}^A$, from~\eqref{f>0} and the fact that $P_r(\infty) = 0$, it follows that  $$s_{n+1}<\infty.$$

From~\eqref{f>0} and the fact that the derivative $G_A'$ can not become arbitrarily large, it follows that $s_k\neq s_j$, $k\neq j$. 

Taking into account~\eqref{G_AIsBounded} and using standard compactness arguments, we can extract a limit $g$-spline $G = \lim\limits_{A\to\infty}G_{A}$, where the convergence of the function and all derivatives of order $\leq r-1$ is uniform on any bounded interval.

Then $G$ has $n+1$ points of oscillation between $-\tau\cdot (f_-+a)$ and $\tau\cdot (f_+-a)$.  Therefore $G$ has $n$ knots, i.~e. $t_k\neq t_j$, $k\neq j$. Really, $G$ has $n$ zeros $z_1<\ldots<z_n$  on $[0,s_{n+1}]$ and $G(\infty) = 0$. Hence $G'$ has $n-1$ zeros on $(0, z_n)$ and one on  $(z_n,\infty)$. Continuing the same way, we obtain that $G^{(r)}$ has $n$ essential sign changes, hence $G$ has $n$ knots. This means that $G\in \Gamma_{n,g}^r[0,\infty)$.

The function $G+\tau a$ is a desired one. Really, it is a perfect $g$-spline with exactly $n$ knots that has $n+1$ points of oscillation between  $-\tau\cdot f_-$ and $\tau\cdot f_+$. Moreover, condition~\eqref{infiniteBoundaryCondition} holds, due to the choice of $a$.
 The lemma is proved.

\begin{lemma}\label{l::derivative'sSignChanges}
Let $n\in\ZZ_+$, $G$ be a perfect $g$-spline from Lemma~\ref{l::halfLineOscillation} with some $\alpha>0$ and $x\in W^r_{\infty, g}[0,\infty)$ be such that 
$$
\|x\|_{C[0,\infty),f_-,f_+} \leq \|G\|_{C[0,\infty),f_-,f_+}.
$$
If $\varepsilon\in (-1,1)$ and $\Delta = G-\varepsilon\cdot x$, then 
$\Delta'$ has at least $n$ sign changes.
\end{lemma}
Let $s_1<\ldots<s_{n+1}$ be the oscillation points of $G$.

We prove the partial case, when $x\equiv 0$, first. Since ${\rm sgn\,}G(s_k) = (-1)^{k+1}$, $k=1,\ldots,n+1$, there exist intervals $I_k\subset (s_k,s_{k+1})$, such that  ${\rm sgn\,}G'(t) = (-1)^k$ for $t\in I_k$, $k=1,\dots, n$. Since $f_\pm$ are non-increasing functions, there exists interval $I_{n+1}\subset (s_{n+1},\infty)$ such that  $G'(t) > 0$ if $n$ is odd or  $G'(t) <0$ if $n$ is even, i.~e.  ${\rm sgn\,}G'(t) = (-1)^{n+1}$, $t\in I_{n+1}$. This implies that $G'$ has at least $n$ sign changes. 

Moreover, by Lemma~\ref{l::singAlternation}, $G'$ has exactly $n$ sign changes and thus 
\begin{equation}\label{g'LastSign}
{\rm sgn\,}G'(t) = (-1)^{n+1},\,t > \inf I_{n+1}.
\end{equation}

Now we return to the proof of the lemma in the general case.

From the conditions of the lemma it follows that $\Delta^{(r)}$ is non-zero almost everywhere and changes sign only in the knots of $G$, hence has exactly $n$ essential sign changes. Then $\Delta'(\infty) = 0$ due to Lemma~\ref{l::derivativeAtInfy2}. Thus $x'(\infty) = 0$, and hence for all $t>0$
\begin{equation}\label{xLessP}
|x'(t)|\leq |P_{r-1}(t)|,
\end{equation}
due to Lemma~\ref{l::GLessP}.

Let $t_n$ be the last knot of $G$. Then $G'(t) = \pm P_{r-1}(t)$, $t>t_n$, since both of them are primitives of either $g$ or $-g$ of order $r-1$ that vanish in infinity together with all of their derivatives. Together with~\eqref{xLessP}, this implies  $|x'(t)|\leq |G'(t)|$, hence 
\begin{equation}\label{lastSignChange}
{\rm sgn\,} \Delta'(t) = {\rm sgn\,} G'(t)
\end{equation}
for   $t>t_n$.

Since ${\rm sgn\,}\Delta(s_k) = (-1)^{k+1}$, $k=1,\ldots,n+1$, there exist intervals $J_k\subset (s_k,s_{k+1})$, such that  ${\rm sgn\,}\Delta'(t) = (-1)^k$ for $t\in J_k$, $k=1,\dots, n$. Combining~\eqref{g'LastSign} and~\eqref{lastSignChange}, we obtain that there exists an interval $J_{n+1}\subset (s_{n+1},\infty)$ such that ${\rm sgn\,}\Delta'(t) = (-1)^{n+1}$ for $t\in J_{n+1}$. Hence $\Delta'$ has at least $n$ sign changes. The lemma is proved.

Applying Lemma~\ref{l::derivative'sSignChanges} with $x\equiv 0$ and Lemma~\ref{l::singAlternation}, we obtain the following lemma.
\begin{lemma}\label{l::rThDerivativeOfGSpline}
If $n\in \ZZ_+$ and $G\in\Gamma^r_{n,g}[0,\infty)$ has $n$ knots and $n+1$ points of oscillation, then ${\rm sgn\,}G^{(r)}(0)= (-1)^r$.
\end{lemma}

\section{Maximally oscillating perfect $g$-spline with given norm}
\subsection{Definition of the functions $\tau_n$}

\begin{lemma}\label{l::partialC_AisWellDefined}
Let $n\in\ZZ_+$ and assume there are two perfect $g$-splines $G_i\in  \Gamma^{r}_{n,g}[0,\infty)$ that oscillate maximally between $-\tau_i\cdot f_-$ and $\tau_i\cdot f_+$, $i=1,2$, and
\begin{equation}\label{ZeroAtInfty}
G_1(\infty) = G_2(\infty),
\end{equation}
then $\tau_1 = \tau_2$.
\end{lemma}

In the case, when $n=0$, $G_1-G_2$ is a constant function, thus due to~\eqref{ZeroAtInfty}, $G_1\equiv G_2$ and the lemma is proved in this case. 

Let $n>0$. Assume the contrary, let for definiteness $\tau_1>\tau_2$. Let $0\leq s_1<\ldots<s_{n+1}$ be the oscillation points of $G_1$, $0<t_1<\ldots<t_n$ be the knots of $G_1$ and $0<u_1<\ldots<u_n$ be the knots of $G_2$. Set $\delta:=G_1-G_2$. Then ${\rm sgn\,} \delta(s_k) = {\rm sgn\,} G_1(s_k)$, $k=1,\ldots, n+1$, and 
\begin{equation}\label{delta(A)}
\delta^{(k)}(\infty) = 0, k=0,\ldots, r-1.
\end{equation}

Hence $\delta$ has $n$ separated zeros on $(0, s_{n+1})$. Due to~\eqref{delta(A)}, $\delta'$ has $n$ separated zeros on $(0, \infty)$. Continuing the same way, we obtain that $\delta^{(r)}$ has $n$ essential sign changes on $(0,\infty)$. This implies, that $\delta^{(r)}$ is not identical zero on each of the intervals $(t_k,t_{k+1})$, $k=0,\ldots, n$, $t_0:=0$, $t_{n+1} := \infty$.

 Hence $u_1<t_1$, since otherwise, due to Lemma~\ref{l::rThDerivativeOfGSpline}, $\delta^{(r)}(t) =0$, $t\in (t_0, t_1)$; $u_2<t_2$, since otherwise $(u_1,u_2)\supset (t_1,t_2)$, and hence $\delta^{(r)}(t) =0$, $t\in (t_1, t_2)$, and so on, $u_n<t_n$. However, we obtain that $\delta^{(r)}(t) = 0$ for $t\in (t_n, t_{n+1})$ and hence $\delta^{(r)}$ has at most $n-1$ sign changes. Contradiction. The lemma is proved.

\begin{lemma}\label{l::C_AisWellDefined}
The conclusion of Lemma~\ref{l::partialC_AisWellDefined} remains true if condition~\eqref{ZeroAtInfty} is substituted by~one of the following conditions:
\begin{enumerate}
\item[(a)] $$\frac{\tau_if_+(\infty) - G_i(\infty)}{\tau_if_-(\infty) + G_i(\infty)} = \alpha,$$
$i=1,2$, $\alpha>0$.
\item[(b)] $n$ is odd and $G_i(\infty) = \tau_if_+(\infty)$, $i=1,2$. 
\item[(c)] $n$ is even and $G_i(\infty) =-\tau_if_-(\infty)$, $i=1,2$. 
\end{enumerate}

\end{lemma}
We prove case (a) first. Set $a := \frac {1}{1+\alpha} f_+(\infty) - \frac {\alpha}{1+\alpha} f_-(\infty)$. Then for $i=1,2$,
$G_i(\infty) =  \tau_i\cdot a$.  The perfect $g$-splines $G_i - G_i(\infty)$ satisfy~\eqref{ZeroAtInfty} and oscillate maximally between $-\tau_i(f_-+a)$ and $\tau_i(f_+ -a)$, $i=1,2$. Lemma~\ref{l::partialC_AisWellDefined} now implies that $\tau_1 = \tau_2$. 

Cases (b) and (c) can be proved similarly. We prove case (b).
If $\tau_1>\tau_2$ and $s_1<\ldots<s_{n+1}$ are the oscillation points of $G_1$, then $G_1(s_{n+1}) = -\tau_1 f_-(s_{n+1}) < G_2(s_{n+1})$ and $G_1(\infty) > G_2(\infty)$. Hence the difference $G_1-G_2$ has $n$ zeroes between the oscillation points of $G_1$ and an additional zero on the interval $(s_{n+1},\infty)$, together at least $n+1$ separated zeros. The contradiction can now be obtained using the arguments from the proof of Lemma~\ref{l::partialC_AisWellDefined}. The lemma is proved.

\begin{remark}
Case (a) of Lemma~\ref{l::C_AisWellDefined} implies the correctness of the following definition.
\end{remark}
\begin{definition}\label{def::C_rn}
For fixed $n\in\ZZ_+$ and $\alpha>0$ set 
$$
	\tau_{n}(\alpha)=\tau_{r,n,g,f_-,f_+}(\alpha):=\|G_\alpha\|_{C_{[0,\infty)},f_-,f_+},
$$
where $G_{\alpha}$ is a spline from $\Gamma^{r}_{n,g}[0,\infty)$ built according to Lemma~\ref{l::halfLineOscillation}.
\end{definition}

\subsection{Some properties of the functions $\tau_n$}

\begin{lemma}\label{l::continuityOfC_n}
For each $n\in\ZZ_+$, $\tau_n$ is continuous on $(0,\infty)$.
\end{lemma}
Let $\beta>0$ and $\beta_s\to\beta$, $s\to\infty$. Denote by $G_\beta$ and $G_s$, $s\in \NN$, the perfect $g$-splines built according to Lemma~\ref{l::halfLineOscillation} with $\alpha$ in the boundary condition~\eqref{infiniteBoundaryCondition} substituted by $\beta$ and $\beta_s$ respectively.

Assume that $\tau_n$ is not continuous at the point $\beta$. Due to Lemma~\ref{l::halfLineOscillation}, the sequence $\left\{\tau_n(\beta_s)\right\}_{s=1}^\infty$ is bounded, hence switching to a subsequence, if needed, we may assume that it has a limit $\lim\limits_{s\to\infty}\tau_n(\beta_s) =\tau\neq \tau_n(\beta)$.
 
Switching to a subsequence, if needed once more, we may assume that the sequences of the knots and the oscillation points of $G_s$ have limits as $s\to\infty$. Moreover, analogously to the proof of Lemma~\ref{l::halfLineOscillation}, we can prove that all oscillation point limits are finite and different. Hence we obtain a perfect $g$-spline $G\in\Gamma_{n,g}^r[0,\infty)$ that oscillates maximally between $-\tau f_-$ and $\tau f_+$, and satisfies~\eqref{infiniteBoundaryCondition} with $\alpha = \beta$. However, this contradicts to statement (a) of Lemma~\ref{l::C_AisWellDefined}. The lemma is proved.

\begin{lemma}\label{l::monotonicityOfC_n}
For $n\in\ZZ_+$ the function $\tau_n$ is non-decreasing in the case of odd $n$, and is non-increasing in the case of even $n$.
\end{lemma}

We prove the lemma for the case of odd $n$; the case of even $n$ can be proved using similar arguments. 

For each $\alpha>0$ denote by $G_\alpha$ a perfect $g$-spline built according to Lemma~\ref{l::halfLineOscillation}. Then, due to~\eqref{infiniteBoundaryCondition},
\begin{equation}\label{valueAtInfy}
G_\alpha(\infty) = \frac{\tau_n(\alpha)}{\alpha+1}\left(f_+(\infty) - \alpha\cdot f_-(\infty)\right).
\end{equation}
Note, that from Lemma~\ref{l::continuityOfC_n} and~\eqref{valueAtInfy} it follows, that $G_\alpha(\infty)$ is a continuous function of $\alpha$.

 Since $n$ is odd, $G_\alpha$ is increasing on the interval $(M_\alpha,\infty)$, where $M_\alpha>0$ is the last point of oscillation of the $g$-spline $G_\alpha$.

Assume the contrary, let $0<\gamma<\beta$ be such that $\tau_n(\gamma) > \tau_n(\beta)$. We claim that 
\begin{equation}\label{inftyPointComparison}
G_\beta(\infty)> G_\gamma(\infty).
\end{equation}
Indeed, otherwise we get an extra zero of $G_\gamma-G_\beta$ on the interval $(M_\gamma,\infty]$ and obtain a contradiction using arguments similar to the proof of Lemma~\ref{l::partialC_AisWellDefined}.

For arbitrary $0<\alpha< \beta$ one has $G_\alpha(\infty)\neq 
G_\beta(\infty)$. Really, otherwise Lemma~\ref{l::partialC_AisWellDefined} implies $\tau_n(\alpha) = \tau_n(\beta)$, and hence, due to~\eqref{valueAtInfy}, $\alpha = \beta$, which is impossible.

From the continuity of $G_\alpha(\infty)$ and~\eqref{inftyPointComparison} it now follows, that 
\begin{equation}\label{Gbeta>Galpha}
G_\beta(\infty) > G_\alpha(\infty)
\end{equation}
for all $0<\alpha < \beta$. 

Note, that equality~\eqref{valueAtInfy} can be rewritten as
$$G_\alpha(\infty) =\tau_n(\alpha)f_+(\infty) -\frac{\alpha}{1+\alpha}\tau_n(\alpha)(f_+(\infty)+f_-(\infty)).$$

This implies that $G_\beta(\infty)<\tau_n(\beta)f_+(\infty)$. Let $\varepsilon> 0$ be such, that $G_\beta(\infty) + \varepsilon<\tau_n(\beta)f_+(\infty)$. Then $\frac{\alpha}{1+\alpha}\tau_n(\alpha)(f_+(\infty)+f_-(\infty))<\varepsilon$ for small enough $\alpha>0$ and hence for such $\alpha$
\begin{gather*}
\tau_n(\alpha)f_+(\infty) = G_\alpha(\infty) + \frac{\alpha}{1+\alpha}\tau_n(\alpha)(f_+(\infty)+f_-(\infty))
\\
< G_\alpha(\infty) + \varepsilon < G_\beta(\infty)+\varepsilon < \tau_n(\beta)f_+(\infty).
\end{gather*}

Thus for all small enough $\alpha>0$ one has  $\tau_n(\alpha) < \tau_n(\beta)$. However, the latter inequality implies that the difference $G_\beta - G_\alpha$ has $n$ separated zeros between oscillation points of $G_\beta$ and additional zero on the interval $(M_\beta,\infty)$ due to~\eqref{Gbeta>Galpha}. Using arguments similar to the proof of Lemma~\ref{l::partialC_AisWellDefined}, we obtain a contradiction. The lemma is proved.

\begin{lemma}\label{l::differentNC_nComparison}
For all $n\in \NN$ and $\alpha,\beta >0$ one has $\tau_n(\alpha)\leq \tau_{n-1}(\beta)$.
\end{lemma}

Assume the contrary, let $\tau_n(\alpha)> \tau_{n-1}(\beta)$. Let $G_\alpha$ and $G_\beta$ be the maximally oscillating splines built according to Lemma~\ref{l::halfLineOscillation} with $n$ and $n-1$ knots respectively. Set $\Delta := G_\alpha - G_\beta$.

Let $s_1<\ldots<s_{n+1}$ be the oscillation points of $G_\alpha$ and let for definiteness $n=2k-1$ be odd, $k\in\NN$. Then for each $m=1,\ldots, 2k-1$, there exists $s_m^1\in (s_m,s_{m+1})$ such that $(-1)^m\Delta'(s_m^1)>0$. Moreover, $G_\alpha$ increases on $(s_{2k},\infty)$ and $G_\beta$ decreases on $(M,\infty)$, where $M$ is the last oscillation point of $G_\beta$. Hence there exists $s_{2k}^1> s_{2k-1}^1$ such that $\Delta'(s_{2k}^1)>0$. Thus $\Delta'$ has $2k-1$ sign changes and moreover, $\Delta'(\infty) = 0$. Hence $\Delta''$ has $2k-1$ sign changes and moreover, $\Delta''(\infty) = 0$. Continuing the same way, we obtain that $\Delta^{(r)}$ has at least $2k-1 = n$ essential sign changes. However, $\Delta^{(r)}$ can change sign only in the knots of $G_\beta$, hence not more than $n-1$ times. Contradiction. The lemma is proved.

\subsection{Asymptotic behavior of the functions $\tau_n$}
\begin{lemma}\label{l::KnotsAndExtremaPositions}
Let $n\in\NN$, $G\in\Gamma_{n,g}^r[0,\infty)$ have $n+1$ points of oscillation between $-f_-$ and $f_+$. Assume $s_1<\ldots <s_n$ are all zeros of $G'$ and $t_1<\ldots<t_n$ are its knots. Then $t_k\geq s_k$, $k=1,\ldots, n$.
\end{lemma}

Assume the contrary, let $k\in \{1,\ldots, n\}$ be such that $t_k<s_k$. Then there exists $\varepsilon>0$ such that $G$ has at most $n-k$ knots on the interval $(s_k-\varepsilon, \infty)$. Since all zeros of $G'$ are simple, on the interval $(s_k-\varepsilon, \infty)$ the function $G'$ changes sign in each of the points $s_l$, $l=k,\ldots, n$, totally $n-k+1$ times. Since $G'(\infty) = 0$, $G''$ has at least $n-k+1$ sign changes on  $(s_k-\varepsilon, \infty)$. Continuing in a similar way, we obtain that $G^{(r)}$ has at least $n-k+1$ sign changes on  $(s_k-\varepsilon, \infty)$. However, this is impossible, since it has only at most $n-k$ knots on this interval. The lemma is proved.

 The following two lemmas describe the limits of the splines from Lemma~\ref{l::halfLineOscillation} as $\alpha\to +0$ and $\alpha\to \infty$. For brevity we write $\tau_n(+0)$ instead of $\lim\limits_{\alpha\to+0}\tau_n(\alpha)$. 

The limits $\tau_n(+0)$, $n\in\NN$, and $\tau_n(\infty) := \lim\limits_{\alpha\to\infty}\tau_n(\alpha)$, $n\in \ZZ_+$, exist, since $\tau_n$ is monotone due to Lemma~\ref{l::monotonicityOfC_n}. Moreover, these limits are positive and finite, due to Lemmas~\ref{l::monotonicityOfC_n} and~\ref{l::differentNC_nComparison}.
\begin{lemma}\label{l::limitsOfG_n}
If $n$ is odd, then there exists a perfect $g$-spline $G_0\in \Gamma^r_{n,g}[0,\infty)$ that has exactly $n$ knots, $n+1$ points of oscillation between $-\tau_n(+0)f_-$ and $\tau_n(+0)f_+$, and such that $G_0(\infty) = \tau_n(+0)f_+(\infty)$. 

If additionally 
\begin{equation}\label{PisBig}
\varliminf\limits_{t\to\infty}\frac {f_-(t)-f_-(\infty)}{|P_r(t)|}  =0,
\end{equation}
then there exists a perfect $g$-spline $G_\infty\in \Gamma_{n-1, g}^r[0,\infty)$ with exactly $n-1$ knots that has $n$ points of oscillation between $-\tau_n(\infty) f_-$ and $\tau_n(\infty) f_+$, and such that 
$$G_\infty(\infty) = -\tau_n(\infty)\cdot f_-(\infty).$$
\end{lemma}
For each $\alpha>0$ denote by $G_\alpha$ a perfect $g$-spline with $n$ knots built according to Lemma~\ref{l::halfLineOscillation}. Let $s_1^\alpha<\ldots< s_{n+1}^\alpha$ be its oscillation points. Since $n$ is odd, $G_\alpha$ is increasing on the interval $(s_{n+1}^\alpha,\infty)$. 

We sketch the proof of the existence of $G_0$. Let $G_0$ be the limiting $g$-spline of appropriately chosen sequence of $G_\alpha$, $\alpha\to+0$, and let $s_{n+1} = \lim\limits_{\alpha\to+0}s_{n+1}^\alpha$. 
Since
$$
G_\alpha(\infty) - G_\alpha(s_{n+1}^\alpha)
\geq G_\alpha(\infty) + \tau_n(\alpha)f_-(\infty),
$$
and the addend $\tau_n(\alpha)f_-(\infty)$ is separated from zero by a constant independent of $\alpha$, then due to Lemma~\ref{l::GLessP}, $s_{n+1}<\infty$. Then, repeating the arguments from Lemma~\ref{l::halfLineOscillation}, one can prove that $G_0$ is the desired $g$-spline.

Next we prove the existence of the $g$-spline $G_\infty$. 

Switching to a subsequence, if needed, let  $s_1\leq\ldots\leq s_{n+1}$ be the (finite or infinite) limits of the sequences $s_1^\alpha,\ldots,  s_{n+1}^\alpha$ as  $\alpha\to \infty$. 
%From~\eqref{f>0} and the fact that the derivative $G_\alpha'$ can not become arbitrarily large, it follows that $s_k\neq s_j$, $k\neq j$. 
Let $G_\infty$ be a limiting  $g$-spline in the sequence $G_\alpha$. 
%Then it has $n+1$ point of oscillation and thus $n$ knots.

Let $z_n^\alpha$ denote the last zero of $G_\alpha'$ and $z_n=\lim\limits_{\alpha\to \infty}z_n^\alpha$ (if needed, we switch to a subsequence once more). We prove that 
\begin{equation}\label{lastOscillationPointLimit}
z_{n}=\infty.
\end{equation}
Assume the contrary, let $z_{n}<\infty$. Then $G_\infty$ increases on $(z_{n},\infty)$, and hence 
$G_\infty(t) = G_\infty(\infty) - |P_r(t)|$ for all big enough $t$, since its restriction to the interval $(M,\infty)$ with enough big $M$ is a primitive of the order $r$ of  either $g$ or $-g$.

From~\eqref{PisBig} it follows, that there exists arbitrarily large  $t$ such that 
$$
\frac {f_-(t)-f_-(\infty)}{|P_r(t)|} <\frac 1 {\tau_n(\infty)}.
$$
%where $\tau:=\|G_\infty\|_{C[0,\infty),f_-,f_+}$.
Hence we obtain 
\begin{gather*}
G_\infty(t) = G_\infty(\infty) - |P_r(t)| = -\tau_n(\infty) f_-(\infty) - |P_r(t)| 
\\< 
\tau_n(\infty) ( -f_-(\infty)-f_-(t)+f_-(\infty)) = - \tau_n(\infty) f_-(t), 
\end{gather*}
which is impossible. Thus~\eqref{lastOscillationPointLimit} holds. Hence $s_{n+1}=\infty$.

Since by Lemma~\ref{l::GLessP} for all $\alpha>0$
\begin{gather*}
|G_\alpha(s^\alpha_{n + 1}) - G_\alpha(s^\alpha_{n})| =  |G_\alpha(s^\alpha_{n + 1}) - G_\alpha(\infty) - (G_\alpha(s^\alpha_{n}) - G_\alpha(\infty))| 
\\ \leq 
|G_\alpha(s^\alpha_{n + 1}) - G_\alpha(\infty)| + |G_\alpha(s^\alpha_{n}) - G_\alpha(\infty)|\leq |P_r(s^\alpha_{n+1})| + |P_r(s^\alpha_n)|,
\end{gather*}
then $s_n^\alpha$ can not become arbitrarily large, due to~\eqref{f>0}. Hence $s_n<\infty$.

From Lemma~\ref{l::KnotsAndExtremaPositions} and~\eqref{lastOscillationPointLimit} it follows, that $G_\infty$ can not have $n$ knots. On the other hand, since $G_\infty$ has $n$ points of oscillation, it can not have less than $n-1$ knots. Hence $G_\infty$ has exactly $n-1$ knots. 

Thus $G_\infty$ is a desired perfect $g$-spline and the lemma is proved.

\begin{lemma}\label{l::limitsOfG_n2}
If $n$ is even, then there exists a perfect $g$-spline $G_\infty\in \Gamma^r_{n,g}[0,\infty)$ that has exactly $n$ knots, $n+1$ points of oscillation between $-\tau_n(\infty)f_-$ and $\tau_n(\infty)f_+$, and such that $G_\infty(\infty) = -\tau_n(\infty)f_-(\infty)$. 

If $n\geq 2$ and additionally 
\begin{equation}\label{pIsBig2}
\varliminf\limits_{t\to\infty}\frac {f_+(t)-f_+(\infty)}{|P_r(t)|}  =0,
\end{equation}
then there exists a perfect $g$-spline $G_0\in \Gamma_{n-1, g}^r[0,\infty)$ with exactly $n-1$ knots that has $n$ points of oscillation between $-\tau_n(+0) f_-$ and $\tau_n(+0) f_+$, and such that 
$$G_0(\infty) = \tau_n(+0)\cdot f_+(\infty).$$

If~\eqref{pIsBig2} holds, then $\tau_0(+0) = \infty$.
\end{lemma}

We prove the last statement. All others can be proved similar to Lemma~\ref{l::limitsOfG_n}.

The spline $G_\alpha$, $\alpha>0$, from Lemma~\ref{l::halfLineOscillation} with $0$ knots has the form $G_\alpha(t) =|P_r(t)| + C_\alpha$, $C_\alpha\in\RR$. If $\tau:=\tau_0(+0) < \infty$, then the splines $G_\alpha$ tend to $G := |P_r| + \tau f_+(\infty)$ for some sequence $\alpha\to+0$. However, condition~\eqref{pIsBig2} implies that for some sufficiently large $t$, $G(t) = |P_r(t)| + \tau f_+(\infty) > \tau f_+(t)$, which is impossible. The lemma is proved.

%\subsection{Relations between the functions $\tau_n$ and $\tau_{n-1}$}
%In Definition~\ref{def::C_rn} we introduced the functions $\tau_n\colon (0,\infty)\to (0,\infty)$ for all $n\in\NN$. Using the function $P_r$ from Definition~\ref{def::Pk}, for convenience we define a (constant) function $\tau_0$.
%\begin{definition}\label{def::C_r0}
%For all $\alpha>0$ set 
%$$
%\tau_0 = \tau_{r,0,g,f_-,f_+}(\alpha):=\inf\limits_{a\in \RR}\|P_r(\cdot) + a\|_{C[0,\infty),f_-,f_+}.
%$$
%\end{definition}
%
%
%
%\begin{lemma}\label{l::limitsOfC_n}
%If 
%\begin{equation}\label{assumption3}
%\varliminf\limits_{t\to\infty}\frac {f_\pm(t)-f_\pm(\infty)}{|P_r(t)|}  =0,
%\end{equation}
%then for all $n\in \NN$, 
%$$\lim\limits_{\alpha\to +0}\tau_{2n-1}(\alpha) = \lim\limits_{\alpha\to +0}\tau_{2n}(\alpha)$$
%and 
%$$\lim\limits_{\alpha\to +\infty}\tau_{2n-1}(\alpha) = \lim\limits_{\alpha\to +\infty}\tau_{2n-2}(\alpha).$$
%\end{lemma}
%The first equality follows from Lemmas~\ref{l::differentNC_nComparison} and~\ref{l::limitsOfG_n2}; the second one follows from Lemmas~\ref{l::differentNC_nComparison} and~\ref{l::limitsOfG_n}. The lemma is proved.

\section{Proof of the main results}
\subsection{Proof of Theorem~\ref{th::ZolotarevSplines}}

We prove statements~(a) and~(b) for the $g$-splines with positively oriented points of oscillation. Due to connection between differently orientated oscillation functions given in Remark~\ref{r::positiveVsNegativeOcsillationOrientation}, the case of negatively oriented points of oscillation can be proved using the same arguments.

Let $n\in\ZZ_+$ be fixed. For $\alpha\in (0,\infty)$ denote by $G_{n,\alpha}$ a perfect $g$-spline with $n$ knots from Lemma~\ref{l::halfLineOscillation}. Set 
$$
T_n:=\left\{\|G_{n,\alpha}\|_{C[0,\infty), f_-,f_+}\colon \alpha\in (0,\infty)\right\}.
$$
 Due to Lemma~\ref{l::continuityOfC_n}, $T_n$ is an interval (open, closed or semi-open).
For all $t\in T_{n+1}$ and $s\in T_n$,  one has $t\leq s$, due to Lemma~\ref{l::differentNC_nComparison}. This implies that $\sup T_{n+1}\leq\inf T_n$. We claim, that $\sup T_{n+1}=\inf T_n$. 

%for each $n\in\ZZ_+$ there exist perfect $g$-splines $G_1, G_2\in \Gamma^r_{n,g}[0,\infty)$ that have $n$ knots and $n+1$ points of oscillation between $-\tau_i f_-$ and $\tau_i f_+$, $i=1,2$, with $\tau_1 = \sup T_{n+1}$ and $\tau_2 = \inf T_n$. Indeed, 

Let $n$ be even for definiteness. Then $n+1$ is odd and we set $G_1$ to be the perfect $g$-spline $G_\infty$ from Lemma~\ref{l::limitsOfG_n} with $n$ substituted by $n+1$; by Lemma~\ref{l::monotonicityOfC_n}, $\tau_{n+1}(\infty) = \sup\limits_{\alpha>0}\tau_{n+1}(\alpha) =\sup T_{n+1}$.

Set $G_2$ to be the  perfect $g$-spline $G_\infty$ from Lemma~\ref{l::limitsOfG_n2}; by Lemma~\ref{l::monotonicityOfC_n}, $\tau_{n}(\infty) = \inf\limits_{\alpha>0}\tau_{n}(\alpha)=\inf T_n$.

Case (c) of Lemma~\ref{l::C_AisWellDefined} now implies that $\sup T_{n+1}=\inf T_n$. 

Since 
$$\bigvee\limits_{0}^\infty G_{n,\alpha} = \int\limits_{0}^\infty|G'_{n,\alpha}(t)|dt,$$ then due to Lemma~\ref{l::GLessP} applied to $G'$, all variations of the splines $G_{n,\alpha}$ are bounded by some number independent of $n$ and $\alpha$. Assumption~\ref{a::fpm} now implies that $\inf T_n\to 0$ as $n\to\infty$.
To finish the proof of statement~(a), it is sufficient to set $a_{n+1}^+:=\inf T_n$, $n\in\ZZ_+$ and note, that $\sup T_0 = \infty$, due to Lemma~\ref{l::limitsOfG_n2}.

To prove statement~(b), it is sufficient to apply Lemma~\ref{l::derivative'sSignChanges} with $x\equiv 0$, Lemma~\ref{l::singAlternation} and note, the $G_\tau^+(0) > 0$ and it decreases on the first monotonicity interval.

To prove statement~(c), we first note that equalities~\eqref{g+AtInfty} and~\eqref{g-AtInfty} hold due to Lemmas~\ref{l::limitsOfG_n} and~\ref{l::limitsOfG_n2}.

Assume the contrary, let $G\in \Gamma^r_{n,g}[0,\infty)$ be such that 
$$
\|G\|_{C[0,\infty), f_-,f_+} < \min\|G_n^\pm\|_{C[0,\infty), f_-,f_+}.
$$
Denote by $G_n$ one of the $g$-splines $G_{a_{n+1}^\pm}^\pm$, such that $${\rm sgn\,} G_n^{(r)}(0) = {\rm sgn\,} G^{(r)}(0).$$ This can be done due to~\eqref{gSplineDerivativesAtZero}.

 Let $s_k$, $k=1,\dots n+1$, be the oscillation points of $G_n$ and $s_{n+2} = \infty$. Set $\delta: = G_n - G$. Then ${\rm sgn\,} \delta(s_k) = {\rm sgn\,} G_n(s_k)$, $k=1,\dots, n+1$. Moreover, this equality also holds for $k=n+2$, due to~\eqref{g+AtInfty} or~\eqref{g-AtInfty}. Using the arguments from the proof of Lemma~\ref{l::partialC_AisWellDefined}, we obtain a contradiction.

 The theorem is proved.
\subsection{Proof of Theorem~\ref{th::InfinityValues}}
The theorem is a combination of the Lemmas~\ref{l::halfLineOscillation}, \ref{l::continuityOfC_n}, \ref{l::monotonicityOfC_n}, \ref{l::differentNC_nComparison}, \ref{l::limitsOfG_n} and~\ref{l::limitsOfG_n2}.
\subsection{Proof of Theorem~\ref{th::modulusOfContinuity}}

Let $x\in W_{f_-,f_+,g}^r[0,\infty)$, $\|x\|_{C[0,\infty),f_-,f_+}\leq \delta$  and assume that 
\begin{equation}\label{contrary}
\|x^{(k)}\|_{C[0,\infty)} >\max\left\{|(G_\delta^+)^{(k)}(0)|,|(G_\delta^-)^{(k)}(0)|\right\}.
\end{equation}

Since $f_\pm$ and $g$ are non-increasing, for arbitrary $\alpha\geq 0$, $x(\cdot +\alpha)\in W_{f_-,f_+,g}^r[0,\infty)$, and 
$\|x\|_{C[0,\infty),f_-,f_+}\geq \|x(\cdot + \alpha)\|_{C[0,\infty),f_-,f_+}$.
Hence without loss of generality we can assume that $\|x^{(k)}\|_{C[0,\infty)} = |x^{(k)}(0)|$. 

Denote by $G$ one of the functions $G_\delta^\pm$, chosen by the condition
\begin{equation}\label{kThDerivative}
{\rm sgn\,} x^{(k)}(0) = {\rm sgn\,} G^{(k)}(0).
\end{equation}
This can be done due to~\eqref{gSplineDerivativesAtZero}. 

Assume for definiteness, that $G= G_\delta^+$, it has $n$ knots and $n+1$ points of oscillation $s_1<\ldots<s_{n+1}$.

There exists $\varepsilon>0$ such that $(1-\varepsilon)|x^{(k)}(0)|>|G^{(k)}(0)|$. Set $\Delta = G- (1-\varepsilon) x$. 

From Lemmas~\ref{l::derivative'sSignChanges} and~\ref{l::singAlternation} it follows, that the functions $\Delta^{(m)}$ and $G^{(m)}$, $m=0,\ldots, r-1$, are $n+1$-piecewise monotone. Moreover, ${\rm sgn\,}\Delta(s_1) = {\rm sgn\,}G(s_1)=1$. Applying the arguments from the proof of Lemma~\ref{l::derivative'sSignChanges}, one can prove that the functions $\Delta$ and $G$ are decreasing on their first interval of monotonicity. Hence from Lemma~\ref{l::singAlternation} it follows, that ${\rm sgn\,}\Delta^{(m)}(0) = {\rm sgn\,}G^{(m)}(0) = (-1)^m$, $m=0,\ldots, r-1$.  However, these equalities with $m=k$ contradict  to~\eqref{contrary} and~\eqref{kThDerivative}.

The theorem is proved.

{\bf Acknowledgments.}
The author would like to thank the referees for their valuable remarks and suggestions.
\bibliographystyle{amsplain}
\bibliography{kovalenko}
\end{document}